\documentclass[11pt]{article}
\usepackage{amsfonts,amscd,amssymb,epsfig}
\newtheorem{definition}{Definition}[section]

\newtheorem{proposition}[definition]{Proposition}
\newtheorem{corollary}[definition]{Corollary}
\newtheorem{remark}[definition]{Remark}

\newtheorem{theorem}[definition]{Theorem}
\newtheorem{example}[definition]{Example}

\def\rawo\lonra{\longrightarrow}

\newenvironment{proof}{{\it Proof.}}{\hfill $ \square $ \vskip 4mm}
\begin{document}
\title{From 3-algebras to  $\Delta$-Groups}
\author
{Mihai D. Staic\thanks{Permanent address: Institute of Mathematics of the  
Romanian Academy, 
PO-Box 1-764, RO-014700 Bucharest, Romania.}\\
Department of Mathematics, SUNY at Buffalo, \\
Amherst, NY 14260-2900, USA\\
e-mail: mdstaic@buffalo.edu}

\date{}
\maketitle
\begin{abstract} We introduce $\Delta$-groups and show how they fit in the context of lattice field theory. To a manifold $M$ we associate a $\Delta$-group $\Gamma(M)$. 
We define the symmetric cohomology $HS^n(G,A)$ of a group $G$ with coefficients in a $G$-module $A$. The $\Delta$-group $\Gamma(M)$ is determined by the action of $\pi_1(M)$ on $\pi_2(M)$ and an element of $HS^3(\pi_1(M),\pi_2(M))$. 
\end{abstract}
\section*{Introduction}
A topological lattice field theory is a prescription of initial data which allows one to construct an invariant for manifolds. 
The machinery works like this: one starts with a triangulation of the manifold, associates a quantity to each simplex, takes an "appropriate" sum and shows that the result is an invariant. Of course, in practice, things are a little bit more complicated.  
The first problem is how you decide, from combinatorial data, if two triangulations  give the same manifold. This can be settled by using an Alexander type theorem \cite{tu} or  Pachner's Theorem \cite{cks}, \cite{ph}. The second problem is to find an algebraic input which reflects the combinatorial equivalence between two triangulations of a manifold. 
In dimension two it is known that topological lattice field theories are in bijection with semisimple associative algebra \cite{fhk}.
For three dimensional manifolds, invariants can be obtained from  various algebraic structures: Hopf algebras in \cite{ku} and \cite{cfs}, 6-j symbols  in \cite{tu} and \cite{fg} or finite groups and 3-cocycles in \cite{dw}. A  review of these results and an example in dimension four can be found in  \cite{cks}.

3-algebras were introduced by Lawrence in \cite{rl} as another possible approach to the problem. A 3-algebra is a vector space $A$ together with three maps $m:A\otimes A\otimes A\to A$, $\overline{m}:A\otimes A\to A\otimes A$ and $P:A\to A$ which satisfy certain compatibilities. Geometrically $m$ represents the projection of  three faces of a tetrahedron to the fourth face, $\overline{m}$ is the projection from two faces of a tetrahedron to the the other two and $P$ in the rotation of a face with an angle of $\frac{2\pi}{3}$. It was shown in \cite{rl} that the 6-j symbol invariant \cite{tu},  fits
naturally in the context of 3-algebras. The first result about 3-algebras is a coherence type theorem  which says that in a 3-algebra a product does not depend on the way we make the evaluation. This is the analog of the fact that in an associative algebra the product does not depend on the way we put the parenthesis.

The idea behind 2-groups can be traced back to \cite{jw} where crossed modules where introduced in connection with the homotopy groups of a topological space.  Formally, a 2-group is a 2-category with one object in which every 1(2)-morphism is invertible. To a topological space $M$ we can associate a 2-group  by taking as 2-morphism maps from the 2-cube $[0,1]\times [0,1]$ to $M$. The horizontal and the vertical composition are obtained by gluing horizontal and 
respectively, vertical two such morphisms. It can be shown that this 2-group is determined by the action of $\pi_1(M)$ on $\pi_2(M)$ and a cohomology class $\alpha \in H^3(\pi_1(M),\pi_2(M))$.  For an up to date description of the subject see \cite{bl}. 
 
In this paper we define strong 3-algebras. These are some particular types of 3-algebras, with $\overline{m}(a\otimes b)=m\otimes id(a\otimes b\otimes u(1))$ where $u:k\to A\otimes A$ is a linear map. The advantage of working with strong 3-algebras is that the relations among $P$, $u$ and $m$ are much simpler. Moreover the examples in \cite{rl} are strong 3-algebras. 

The set theoretical equivalent for 3-algebras is the notion of 
a $\Delta$-group. 
We give a construction which associates to a manifold $M$ a $\Delta$-group $\Gamma(M)$. This  generalization is in the spirit of the definition for the fundamental group $\pi_1(M)$. The idea is to replace paths between based points with  2-paths between "based curves" (in other words, equivalences classes of maps from a 2-simplex to $M$ which restricted to boundary are certain fixed curves). We show that every finite $\Delta$-group give rise to a strong 3-algebra.
We study a certain class of $\Delta$-groups  associated to a $G$ module $A$ and show that they are classified by the symmetric cohomology $HS^3(G, A)$.
In particular  we associate to every manifold an element in 
$HS^3(\pi_1(M), \pi_2(M))$. 

In the last section we  define the symmetric cohomology $HS^n(G,A)$. For this we give an action of the symmetric group $\Sigma_{n+1}$ on $C^n(G,A)$ and show that it is compatible with the usual differential. This is similar with way one defines  the cyclic cohomology for algebras (of course in that case one uses the action of the cyclic group $C_{n+1}$). 

The problems we study in this paper can be broadly classified as belonging to the field of ``2-algebra''. Namely we look to maps between morphisms. Our approach is different then the classical one in looking to 2-simplexes rather then 2-cubes.  This makes it more natural to talk about tri-products than vertical and horizontal compositions. 
It also allows us to point out a certain symmetry that is overlooked  by the definition of 2-groups.

\section{Preliminaries}

In this section we recall a few definitions  and results about 3-algebras. For more details we refer to \cite{rl}. In what follows $k$ is a field, $\otimes$ means  $\otimes_k$. If $V$ is a vector space $\tau_{i,j}:V^{\otimes n}\to V^{\otimes n}$ is the the transposition which interchanges the $i$-th and  $j$-th positions. 
\begin{definition} A 3-algebra over $k$ is a vector space $A$ endowed with 
$k$-linear maps, 
$$P:A\to A \; \; \; \; \; {\rm( of \; order\;  3, \; P^3=id)}$$
$$m:A \otimes  A\otimes A\to A$$
$$\overline{m} :A\otimes A\to A\otimes A$$
which satisfy the following conditions:

${\rm(i)}\;  m(m\otimes 1\otimes 1)=m(1\otimes 1\otimes m)\tau_{34}(1\otimes 
\overline{m} \otimes 1 \otimes 1)\tau_{34}$

${\rm (ii)}\; (1\otimes m)\tau_{23}(\overline{m}\otimes 1\otimes 1)=
\overline{m}(1 \otimes m)\tau_{12}(P^{-1}\otimes 1\otimes 1\otimes 1)(\overline{m} \otimes 1\otimes 1)(P\otimes P\otimes 1\otimes 1)\tau_{23}$

${\rm (iii)}\; \overline{m}(m\otimes 1)= (1\otimes m)\tau_{12}(P^2\otimes \overline{m} \otimes 1)(1\otimes 1 \otimes \overline{m})\tau_{12}\tau_{23}$

${\rm (iv)}\; (1\otimes \overline{m})\tau_{12}(1\otimes \overline{m})=
(\overline{m}\otimes 1)(1\otimes \overline{m})(P\otimes P\otimes 1)(\overline{m} \otimes 1)(1\otimes P^{-1}\otimes 1)$

${\rm (v)}\; (1\otimes m)\tau_{23}(\overline{m} \otimes P^2 \otimes 1)=
(m\otimes 1)(1\otimes 1\otimes \overline{m})$

${\rm (vi)} \; Pm=m(P\otimes P\otimes P)\tau_{23}\tau_{12}$

${\rm (vii)}\; \overline{m} \; commutes \; with \;  (P^2\otimes P)\tau_{12}$
\end{definition}

\begin{definition}
A 3-algebra is said to be orthogonal if:

${\rm (viii)} \; (1\otimes P^2)\tau_{12}\overline{m}(P\otimes P)\overline{m}= Q:A\otimes A\to A\otimes A$ is a projection and $m$ vanishes on 
$(ker Q)\otimes A$. 
\end{definition}

\begin{figure}[htbp] 
\hspace*{-80pt}
\centerline{\psfig{file=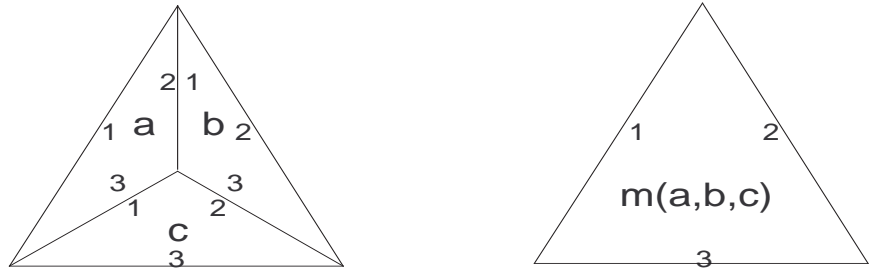,height=3.5cm,width=4.5cm}} 
\caption{$m(a\otimes b\otimes c)$}
\label{fig2}
\end{figure}

\begin{figure}[htbp] 
\begin{center}
\centerline{\psfig{file=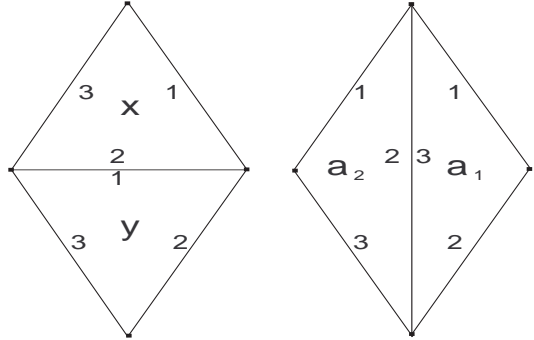,height=4.5cm,width=7cm}} 
\caption{$\overline{m}(x\otimes y)=\sum a_1\otimes a_2$}
\label{fig22}
\end{center}
\end{figure}

The geometric pictures for $m$ and $\overline{m}$ are figures \ref{fig2} and 
\ref{fig22} respectively.

A product in a 3-algebra is a labeled triangulation $\Pi$ of a triangle $T$.  More exactly, each triangle of the triangulation has labels $1, 2, 3$ placed on sides  and is labeled with an element of $A$.   
An ordered evaluation ${\cal T}$ of a  product $\Pi$ is  sequence of triangulations, starting with $\Pi$ and ending with the trivial triangulation of $T$ such that, at each step we change the triangulation by a move depicted in Figure \ref{fig2} or Figure \ref{fig22}.

\begin{theorem} \cite{rl}
Suppose that $A$ is an orthogonal 3-algebra and $\Pi$ is labeled triangulation of a triangle. Then the composition of $m$, $\overline{m}$ and $P$ specified 
by an evaluation ${\cal T}$ of $\Pi$, has a image in $A$ which is independent of the choice of ${\cal T}$.
\end{theorem}

Let $I$ be a set, $f:I^6\to k$ a map which will be denoted 
$$(a,b,c,i,j,k)\to 
\left\vert\begin{array}{ccc}
a&b&c\\
i&j&k
\end{array}\right\vert
$$
We assume that $f$ is invariant under the action of the symmetric group $\Sigma_4$. Also let $w:I\to k$.
\begin{example} \cite{rl} The $k$ linear space A, generated by $\{ e_{ijk}\vert i,j,k\in I\}$ together with maps:
$$
P(e_{ijk})=e_{jki}
$$
$$
m(e_{akj}\otimes e_{k'bi}\otimes e_{j'i'c})=\delta_{ii'}\delta_{jj'}\delta_{kk'}
\left\vert\begin{array}{ccc}
a&b&c\\
i&j&k
\end{array}\right\vert
e_{abc}
$$
$$\overline m(e_{j_2bc}\otimes e_{b'aj_1})=\delta_{bb'}\sum_jw_j^2
\left\vert\begin{array}{ccc}
j_2&a&j\\
j_1&c&b
\end{array}\right\vert
e_{j_2aj}\otimes e_{cjj_1}
$$
defines a 3-algebra if and only if for all $a,b,c,e,f,j_1,j_2,j_3,j_{23}\in I$ we have:
$$
\sum_jw_j^2
\left\vert\begin{array}{ccc}
e&j_3&j\\
j_2&a&j_{23}
\end{array}\right\vert
\left\vert\begin{array}{ccc}
j&c&j_1\\
b&a&j_2
\end{array}\right\vert
\left\vert\begin{array}{ccc}
j_3&c&f\\
j_1&e&j
\end{array}\right\vert
=\left\vert\begin{array}{ccc}
j_3&c&f\\
b&j_{23}&j_2
\end{array}\right\vert
\left\vert\begin{array}{ccc}
e&f&j_1\\
b&a&j_{23}
\end{array}\right\vert
$$
\label{example1}
\end{example}

\section{Strong 3-Algebras}

In this section we study a particular type of 3-algebras for which there is a stronger relation between $m$ and $\overline{m}$. The geometric interpretation of this dependence is depicted in Figure \ref{fig3}.

\begin{figure}[htbp] 
\hspace*{-50pt}
\centerline{\psfig{file=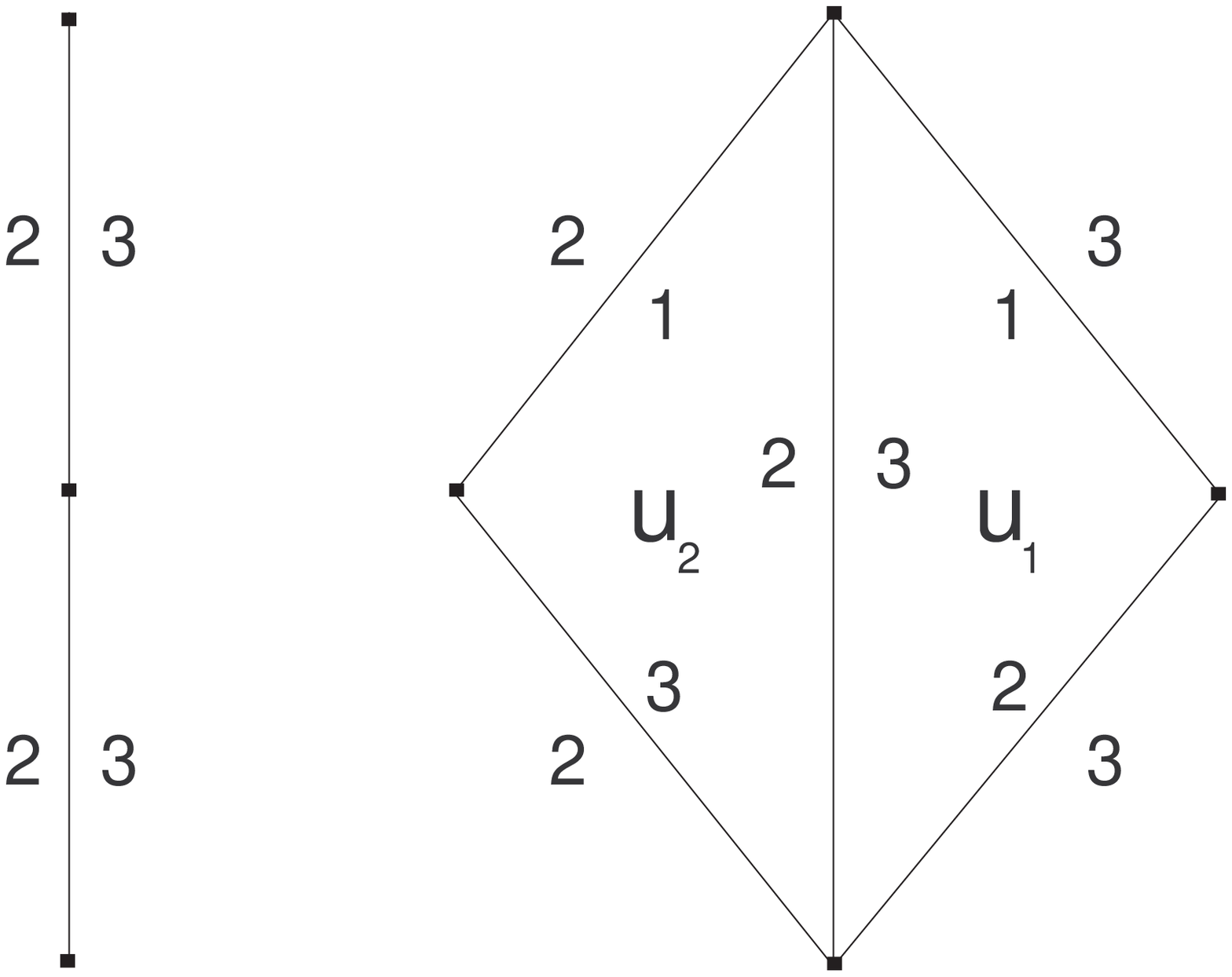,height=4.5cm,width=5.5cm}} 
\caption{$u(1)=\sum u_1 \otimes u_2=\sum P^2(u_2)\otimes P(u_1)$}
\label{fig1}
\end{figure}

\begin{figure}[htbp] 
\hspace*{-40pt}
\centerline{\psfig{file=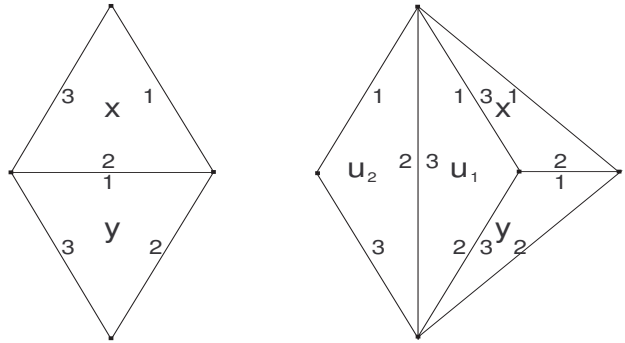,height=4.5cm,width=7cm}} 
\caption{$\tilde m(x\otimes y)=\sum m(x\otimes y\otimes u_1)\otimes u_2$}
\label{fig3}
\end{figure}

\begin{definition}
A strong 3-algebras over a field $k$ is a  vector space $A$ with $k$-linear maps:
$$P:A\to A$$
$$u:k\to A\otimes A$$
$$m:A\otimes A\otimes A \to A$$ 
such that $(A, P, m, \tilde m)$ is a 3-algebra, where $\tilde m : A\otimes A\to A\otimes A$ is defined by 
$$\tilde m(a\otimes b)=(m\otimes id)(a\otimes b \otimes u(1))$$
\end{definition}
Figure \ref{fig1} suggests that $u$ should have a certain symmetry.  We assume that  we have the following identity:
\begin{equation}
\sum u_1 \otimes u_2=\sum P^2(u_2)\otimes P(u_1)\label{s2}
\end{equation}

\begin{figure}[htbp] 
\hspace*{-80pt}
\centerline{\psfig{file=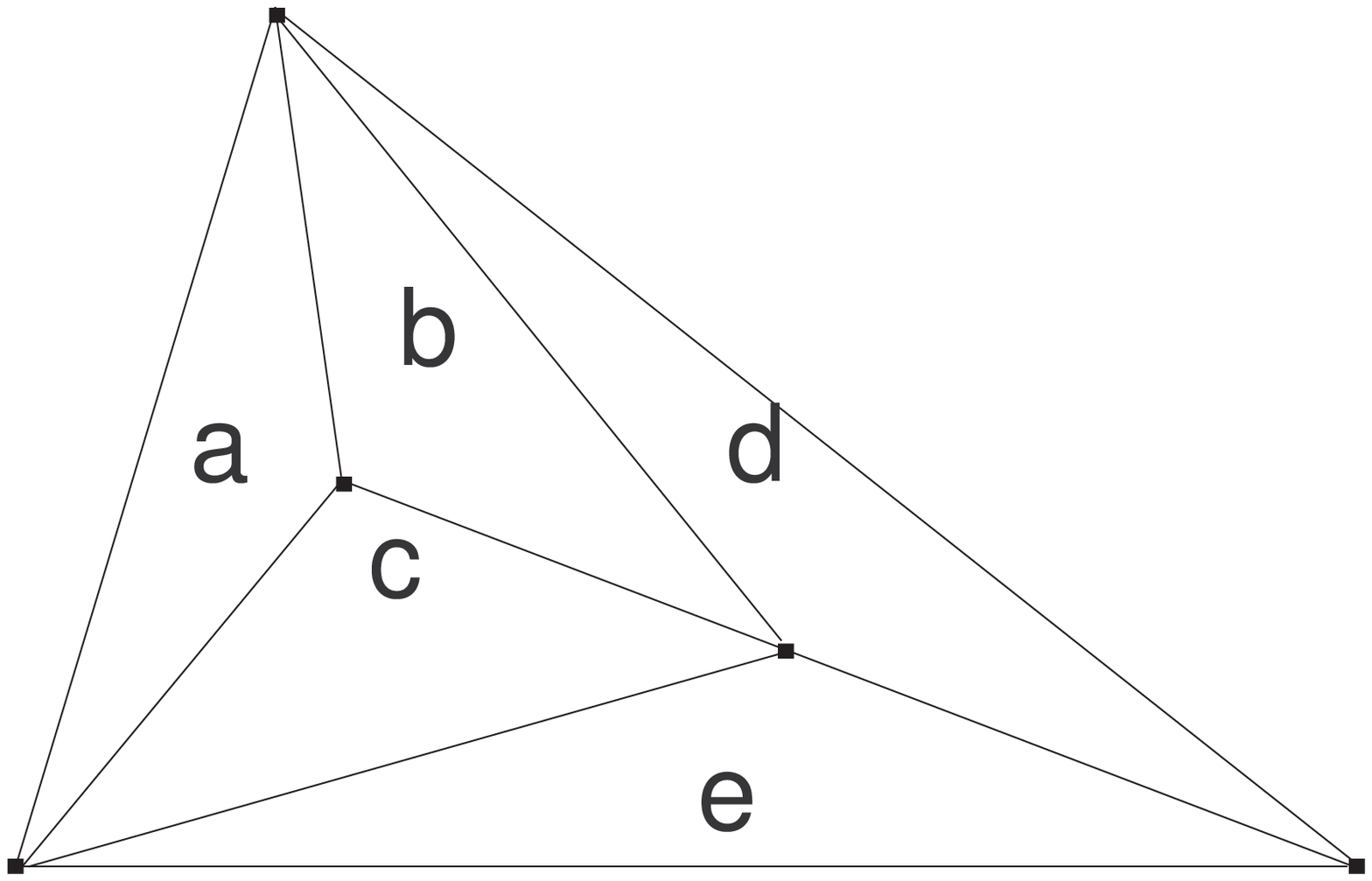,height=3cm,width=4.5cm}} 
\caption{$m(m(a\otimes b\otimes c)\otimes d\otimes e)=\sum m(a\otimes m(b \otimes d\otimes u_1)\otimes m(c\otimes u_2 \otimes e))$}
\label{fig4}
\end{figure}

\begin{figure}[htbp] 
\hspace*{-80pt}
\centerline{\psfig{file=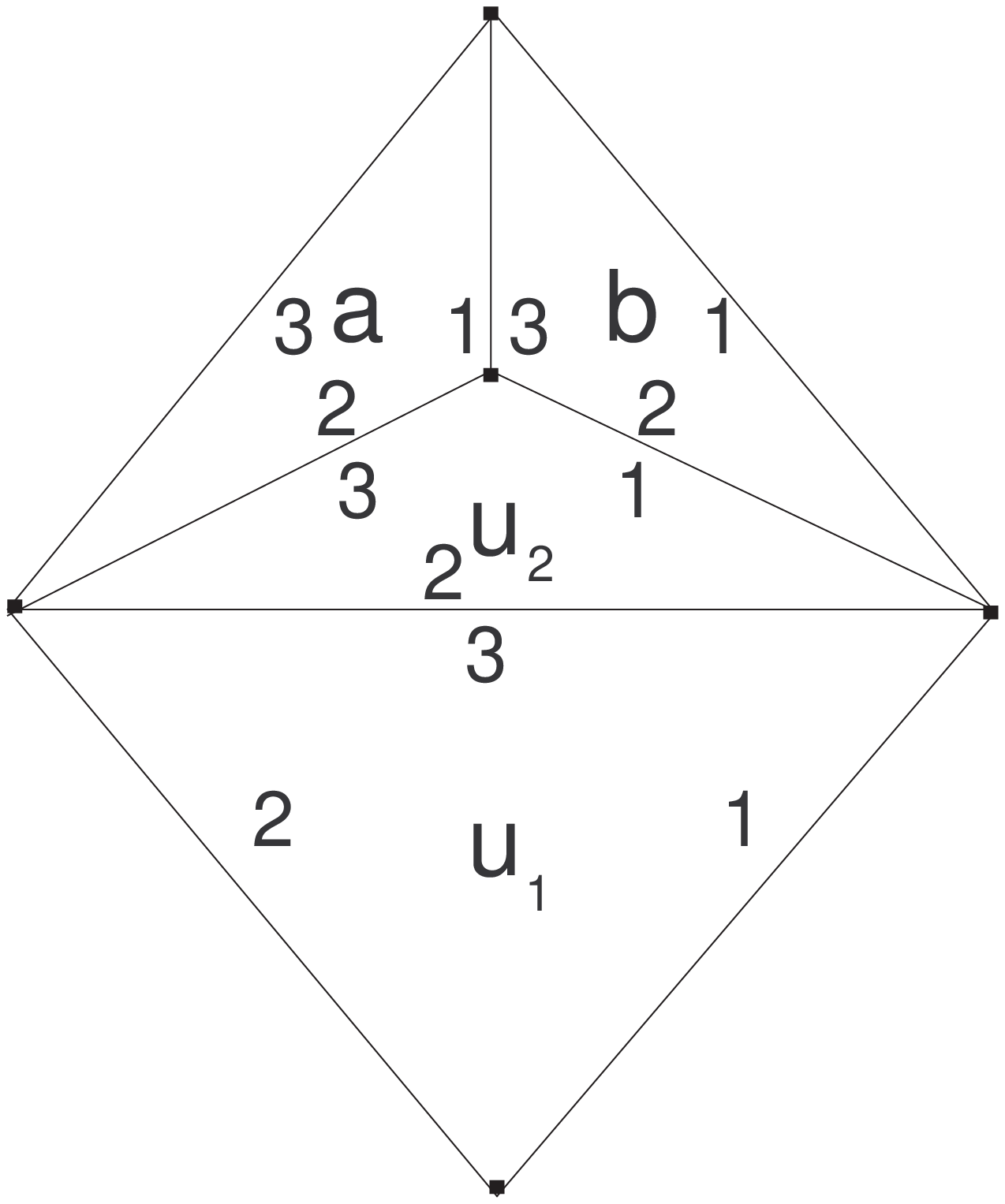,height=5cm,width=4.2cm}} 
\caption{$\sum u_1 \otimes m(b\otimes u_2 \otimes a)=m(P(b)\otimes a\otimes u_1)\otimes u_2$}
\label{fig5}
\end{figure}

It is natural to ask what are necessary and sufficient conditions which make $(A,m,u,P)$ a strong 3-algebra on $A$. The answer is given in the next proposition.
\begin{proposition}
Let $A$ be a vector space over $k$ and let $P:A\to A$, $u:k\to A\otimes A$ and $m:A\otimes A\otimes A \to A$ be $k$-linear maps. Suppose that $u$ satisfies \ref{s2}, then $(A, P, u, m)$ is a strong 3-algebra if and only if:
\begin{equation}
P^3=id\label{s1}
\end{equation}
\begin{equation}
P(m(a\otimes b\otimes c))=m(P(b)\otimes P(c)\otimes P(a))\label{s3}
\end{equation}
\begin{equation}
\sum m(P(b)\otimes a\otimes u_1)\otimes u_2= \sum u_1 \otimes m(b\otimes u_2 \otimes a) \label{s4}
\end{equation}
\begin{equation}
m(m(a\otimes b \otimes c)\otimes d\otimes e)=
\sum m(a\otimes m(b \otimes d\otimes u_1)\otimes m(c\otimes u_2 \otimes e))\label{s5}
\end{equation}
\end{proposition}

\begin{proof} Since our formulas involve several copies of $u(1)$ in the same time we shall use the following notations $u(1)=\sum u_1\otimes u_2=\sum U_1\otimes U_2=\sum \tilde u_1 
\otimes \tilde u_2$.

Obviously (vi) and \ref{s3} are the same. We have:
\begin{eqnarray*}
& &m(1\otimes 1\otimes m)\tau_{34}(1\otimes \tilde m \otimes 1 \otimes 1)\tau_{34}(a\otimes b\otimes c\otimes d\otimes e)=\\
&&=\sum m(1\otimes 1\otimes m)\tau_{34}(a\otimes m(b\otimes d\otimes u_1)\otimes u_2 \otimes c \otimes e)\\ 
&&=\sum m(a\otimes m(b \otimes d\otimes u_1)\otimes m(c\otimes u_2 \otimes e))
\end{eqnarray*} 
And so (i) follows from \ref{s5}. We compute:
\begin{eqnarray*}
\tilde m(P^2\otimes P)\tau_{12}(a\otimes b)&=&\tilde m(P^2(b) \otimes P(a))\\
&=&\sum m(P^2(b)\otimes P(a)\otimes u_1)\otimes u_2
\end{eqnarray*}
\begin{eqnarray*}
(P^2 \otimes P)\tau_{12}\tilde m(a\otimes b)&=&\sum (P^2 \otimes P)\tau_{12}(m(a\otimes b \otimes u_1)\otimes u_2)\\
&=&\sum P^2(u_2)\otimes P(m(a\otimes b\otimes u_1))\\
&=&\sum P^2(u_2)\otimes m(P(b) \otimes P(u_1)\otimes P(a))\\
&=&\sum u_1\otimes m(P(b) \otimes u_2 \otimes P(a))
\end{eqnarray*}
and now (vii) follows from \ref{s4}.  To prove (v) we use \ref{s4} and the following two equations
\begin{eqnarray*}
(1\otimes m)\tau_{23}(\tilde m \otimes P^2 \otimes 1)(a\otimes b\otimes c\otimes d)&=&\sum (1\otimes m)\tau_{23}(m(a\otimes b\otimes u_1)\otimes u_2\otimes P^2(c) \otimes d)\\
&=&\sum m(a\otimes b\otimes u_1)\otimes m(P^2(c)\otimes u_2\otimes d)
\end{eqnarray*}
\begin{eqnarray*}
(m\otimes 1)(1\otimes 1\otimes \tilde m)(a\otimes b\otimes c\otimes d)&=&\sum (m\otimes 1)(a \otimes b\otimes m(c\otimes d \otimes u_1)\otimes u_2)\\
&=&\sum m(a\otimes b\otimes m(c\otimes d\otimes u_1))\otimes u_2
\end{eqnarray*}

We compute:
\begin{eqnarray*}
(1\otimes \tilde m)\tau_{12}(1\otimes \tilde m)(a \otimes b\otimes c)&=&\sum (1\otimes \tilde m)(m(b\otimes c\otimes u_1) \otimes a\otimes u_2)\\
&=&\sum m(b\otimes c\otimes u_1)\otimes m(a\otimes u_2 \otimes U_1)\otimes U_2
\end{eqnarray*}
\begin{eqnarray*}
&&(\tilde m \otimes 1)(1\otimes \tilde m)(P\otimes P \otimes 1)(\tilde m\otimes 1)(a \otimes P^{-1}(b) \otimes c)=\\
&&=\sum (\tilde m \otimes 1)(1\otimes \tilde m)(P\otimes P \otimes 1)(m(a\otimes P^2(b) \otimes u_1)\otimes u_2 \otimes c) \\
&&=\sum (\tilde m \otimes 1)(1\otimes \tilde m)(m(b\otimes P(u_1)\otimes P(a))\otimes P(u_2) \otimes c)\\
&&=\sum (\tilde m \otimes 1)(m(b\otimes P(u_1)\otimes P(a))\otimes m(P(u_2)\otimes c \otimes U_1)\otimes U_2)\\
&&=\sum m(m(b\otimes P(u_1)\otimes P(a))\otimes m(P(u_2)\otimes c \otimes U_1)\otimes \tilde u_1) \otimes \tilde u_2 \otimes U_2\\
&&=\sum m(b\otimes c\otimes m(P(a)\otimes U_1\otimes \tilde u_1))\otimes \tilde u_2 \otimes U_2\\
&&=\sum m(b\otimes c\otimes u_1)\otimes m(a\otimes u_2\otimes U_1) \otimes U_2
\end{eqnarray*}
so we get (iv). For (iii) we check the following  two equalities:
\begin{eqnarray*}
\tilde m(m\otimes 1(a\otimes b\otimes c \otimes d))&=&\tilde m(m(a\otimes b\otimes c)\otimes d)\\
&=&\sum m(m(a\otimes b\otimes c)\otimes d\otimes u_1)\otimes u_2\\
&=&\sum m(a\otimes m(b\otimes d\otimes U_1)\otimes m(c\otimes U_2\otimes u_1))
\otimes u_2\\
&=&\sum m(a \otimes m(b\otimes d\otimes U_1)\otimes u_1)\otimes m(P^2(c)\otimes 
u_2\otimes U_2) 
\end{eqnarray*}
\begin{eqnarray*}
&&(1\otimes m)\tau_{12}(P^2 \otimes \tilde m \otimes 1)(1\otimes 1\otimes 
\tilde m)\tau_{12}\tau_{23}(a\otimes b\otimes c\otimes d)=\\
&&=\sum 1\otimes m(\tau_{12}(P^2 \otimes \tilde m \otimes 1)(1\otimes 1\otimes 
\tilde m)(c\otimes a\otimes b\otimes d))\\
&&=\sum 1\otimes m(\tau_{12}(P^2 \otimes \tilde m \otimes 1)(c\otimes a\otimes 
m(b\otimes d\otimes u_1)\otimes u_2)\\
&&=\sum 1\otimes m(\tau_{12}(P^2(c)\otimes m(a\otimes m(b\otimes d\otimes u_1)\otimes U_1)\otimes U_2 \otimes u_2)\\
&&=\sum m(a\otimes m(b\otimes d\otimes u_1) \otimes U_1)\otimes m(P^2(c)\otimes 
U_2\otimes u_2)
\end{eqnarray*}
Finally we have:
\begin{eqnarray*}
&&(1\otimes m)\tau_{23}(\tilde m \otimes 1\otimes 1)(a\otimes b\otimes c\otimes d)=\\
&&=\sum (1 \otimes m)\tau_{23}(m(a\otimes b\otimes u_1)\otimes u_2\otimes c\otimes d)\\
&&=\sum (1\otimes m)(m(a\otimes b\otimes u_1)\otimes c\otimes u_2\otimes d)\\
&&=\sum m(a\otimes b\otimes u_1)\otimes m(c\otimes u_2\otimes d)\\
&&=\sum m(a\otimes b\otimes m(P(c)\otimes d\otimes u_1))\otimes u_2
\end{eqnarray*}
and
\begin{eqnarray*}
&&\tilde m(1\otimes m)\tau_{12}(P^2\otimes 1\otimes 1\otimes 1)(\tilde m\otimes 1\otimes 1)(P\otimes P\otimes 1\otimes 1)\tau_{23}(a\otimes b\otimes c\otimes d)=\\
&&=\tilde m(1\otimes m)\tau_{12}(P^2\otimes 1\otimes 1\otimes 1)(\tilde m\otimes 1\otimes 1)(P(a)\otimes P(c)\otimes b\otimes d)\\
&&=\sum \tilde m(1\otimes m)\tau_{12}(P^2\otimes 1\otimes 1\otimes 1)(m(P(a)\otimes P(c)\otimes u_1)\otimes u_2\otimes b\otimes d)\\
&&=\sum \tilde m(1\otimes m)\tau_{12}(m(P^2(u_1)\otimes a\otimes c)\otimes u_2\otimes b\otimes d)\\
&&=\sum \tilde m(1\otimes m)(u_2\otimes m(P^2(u_1)\otimes a\otimes c)\otimes b\otimes d)\\
&&=\sum \tilde m(u_2\otimes m (m(P^2(u_1)\otimes a\otimes c)\otimes b\otimes d))\\
&&=\sum m(u_2\otimes m (m(P^2(u_1)\otimes a\otimes c)\otimes b\otimes d)\otimes U_1)\otimes U_2\\
&&=\sum m(m(a\otimes u_2\otimes P(c)) \otimes m (P^2(u_1)\otimes b\otimes d)\otimes U_1)\otimes U_2\\
&&=\sum m(a\otimes b\otimes m(P(c)\otimes d\otimes U_1))\otimes U_2
\end{eqnarray*}
Which completes our proof.
\end{proof}

\begin{example} Consider the construction from Example \ref{example1}. We define $u:k\to A\otimes A$ 
$$u(1)=\sum_{j,u,v} w_j^2e_{uvj}\otimes e_{ujv}$$
Then $(A,m,u,P)$ is a strong 3-algebra.
\end{example}

The following example is inspired by the Dijkgraaf-Witten invariant associated to a finite group $G$ and a 3-cocycle \cite{dw}.
\begin{example}
Let $G$ be a finite group and $\alpha:G\times G\times G\to k$ a 3-cocycle. We consider the vector space $k[G^{(3-1)}]$ which  has a basis indexed by the the triples $(g,h,k)$ with the property $khg=1$ (notice that $k[G^{(3-1)}]$ has dimension $\vert G\vert ^2$). 

Define three linear maps $P:k[G^{(3-1)}]\to k[G^{(3-1)}]$, $u:k \to k[G^{(3-1)}]\otimes k[G^{(3-1)}]$ and $m:k[G^{(3-1)}]\otimes k[G^{(3-1)}]\otimes k[G^{(3-1)}]\to k[G^{(3-1)}]$  determined by:
$$P((g,h,k))=(h,k,g)$$
$$u(1)=\sum_{g,h}(g,h, (hg)^{-1})\otimes (g^{-1}, hg, h^{-1})$$
$$m((x,y,z),(p,q,r),(a,b,c))=\delta(az)\delta(br)\delta(py)\alpha(z,r,q)
(x,q,c)$$
Then $k[G^{(3-1)}]$ is a strong 3-algebras if and only if we have: $\alpha(g,h,k)=\alpha(gh,k,(hk)^{-1})=\alpha((hk)^{-1},g^{-1}, gh)=\alpha(hk,k^{-1},(gh)^{-1})$.
\label{dwe}
\end{example}

\section{$\Delta$-Groups and Manifolds}

It is well known that to every finite group  $G$ one can associates the group algebra $kG$. It would be nice to have a similar construction for 3-algebras. For this we need to replace groups with some other set theoretical structure.
 
In this section we associate to every manifold $M$ a $\Delta$-group $\Gamma(M)$ which will be, in some sense, the higher dimensional analog of the fundamental group $\pi_1(M)$. We first describe the construction and then give the actual definition for $\Delta$-groups. We conclude by showing that from every finite $\Delta$-group one can construct a strong 3-algebra.

Let $M$ be a manifold such that no element from $\pi_1(M)$ 
has order 2. Take $m_0$ a base point in  $M$. Let $\Omega(M,m_0)$ the set of all closed paths starting at $m_0$. Consider the map $pr:\Omega(M,m_0)\to \pi_1(M)$  that sends a path to its homotopy class. We fix a section $s:\pi_1(M)\to \Omega(M,m_0)$ satisfying these two conditions:
\begin{eqnarray*}
s(\alpha^{-1})=s(\alpha)\circ(t\to 1-t)\\
s(1)={\rm constant}\; {\rm map} \; (t\to m_0)
\end{eqnarray*}
Set ${\cal B}(M)=s(\pi_1(M))$. 

Consider the  standard 2-simplex $\Delta_2=\{(x_0,x_1,x_2)\vert x_0+x_1+x_2=1, x_i\ge 0\}$. We denote by $[0]$, $[1]$ and $[2]$ the three vertexes of the simplex and by $[0,1]$, $[1,2]$ and $[0,2]$ the corresponding edges. 

For $\alpha$, $\beta \in {\cal B}(M)$ we define $\Gamma(\alpha, \beta)$ to be the set of homotopy equivalence classes of maps $a:\Delta_2 \to M$ such that $a_{\vert [1,0]}=\alpha$, $a_{\vert [2,1]}=\beta$ and 
$a_{\vert [2,0]}=\beta\alpha$. Here by $\beta\alpha$  we mean the element $s(pr(\beta)pr(\alpha))\in{\cal B}(M)$. $\Gamma(\alpha, \beta)$ is never empty because  in $\pi_1(M)$ we have     $pr(\beta\alpha)^{-1}pr(\beta)pr(\alpha)=1$.

Let $p,q: \Delta_2 \to \Delta_2$ be the maps defined by $$p(x_0,x_1,x_2)=(x_1,x_2,x_0) \; {\rm  and} \;  q(x_0,x_1,x_2)=(x_1,x_0,x_2)$$ Notice that $p^3=id_{\Delta}$, $q^2=id_{\Delta}$ and that they give a representation of the symmetric group $\Sigma_3$. 

We define $P:\Gamma(\alpha, \beta)\to \Gamma(\beta, (\beta\alpha)^{-1})$, by  $P(a)=ap$ and $Q:\Gamma(\alpha, \beta)
\to \Gamma(\alpha^{-1},\beta\alpha)$, by $Q(a)=aq$. 
Obviously we have: 
$$P^3=id, Q^2=id \; {\rm and} \;  QP=P^2Q$$
Consider now the $3$-dimensional simplex $\Delta_3=\{(x_0,x_1,x_2,x_3)\vert x_0+x_1+x_2+x_3=1, x_i\ge 0\}$. Take $a\in \Gamma(\alpha,\beta^{-1})$, $b\in \Gamma(\beta,\gamma^{-1})$ and $c\in \Gamma(\beta^{-1}\alpha, \gamma^{-1}\beta)$
 
We want to define a map $\omega:\Delta_3\to M$ by gluing $a$, $b$ and $c$ on three faces and then extending the map to the rest of the simplex. First we define  $\omega_{\vert [1,2,0]}=a$.  Since $a\in \Gamma(\alpha, \beta^{-1})$ it means that $\omega_{\vert [0,2]}=\beta^{-1}$ and because $b\in \Gamma(\beta,\gamma^{-1})$ and $\beta^{-1}(t)=\beta(1-t)$ we can extend 
$\omega$ such that $\omega_{\vert [0,2,3]}=b$. Using a similar argument we can assume that $\omega_{\vert [1,0,3]}=c$. Finally, we extend  $\omega $ to the whole $\Delta_3$.  From construction we can see that homotopy class of $\omega_{\vert [1,2,3]}$ depends only on the homotopy class of $a$, $b$ and $c$. 
It means that we have defined a map 
$$m:\Gamma(\alpha,\beta^{-1})\times \Gamma(\beta,\gamma^{-1})\times \Gamma(\beta^{-1}\alpha, \gamma^{-1}\beta)\to \Gamma(\alpha, \gamma^{-1})$$
$$m(a,b,c)=\omega_{\vert[1,2,3]}$$

For every  $\alpha$, $\beta\in {\cal B}(M)$ and $f\in \Gamma(\alpha,\beta)$  we construct $U(f)\in \Gamma(\alpha^{-1}, \beta\alpha)$ in the following way: first we define $\theta: \Delta_3 \to M$ such that $\theta_{\vert [1,2]}=m_0$. For every point $x$ on the edge $[1,2]$ we put 
$\theta_{\vert [x,0]}=\alpha$ and  $\theta_{\vert [3,x]}=\beta$. We extend $\theta$ such that $\theta_{\vert [0,1,3]}=f$.  Up to homotopy any extension of $\theta$ to the whole $\Delta_3$ will give the same restriction  on $[2,0,3]$. 
And so we have a map: 
$$U: \Gamma(\alpha,\beta) \in \Gamma(\alpha^{-1}, \beta\alpha)$$ 
$$U(f)=\theta_{\vert [2,0,3]}$$
It is not difficult to see that $U(f)=Q(f)$.

\begin{figure}[htbp] 
\vspace*{-83pt}
\centerline{\psfig{file=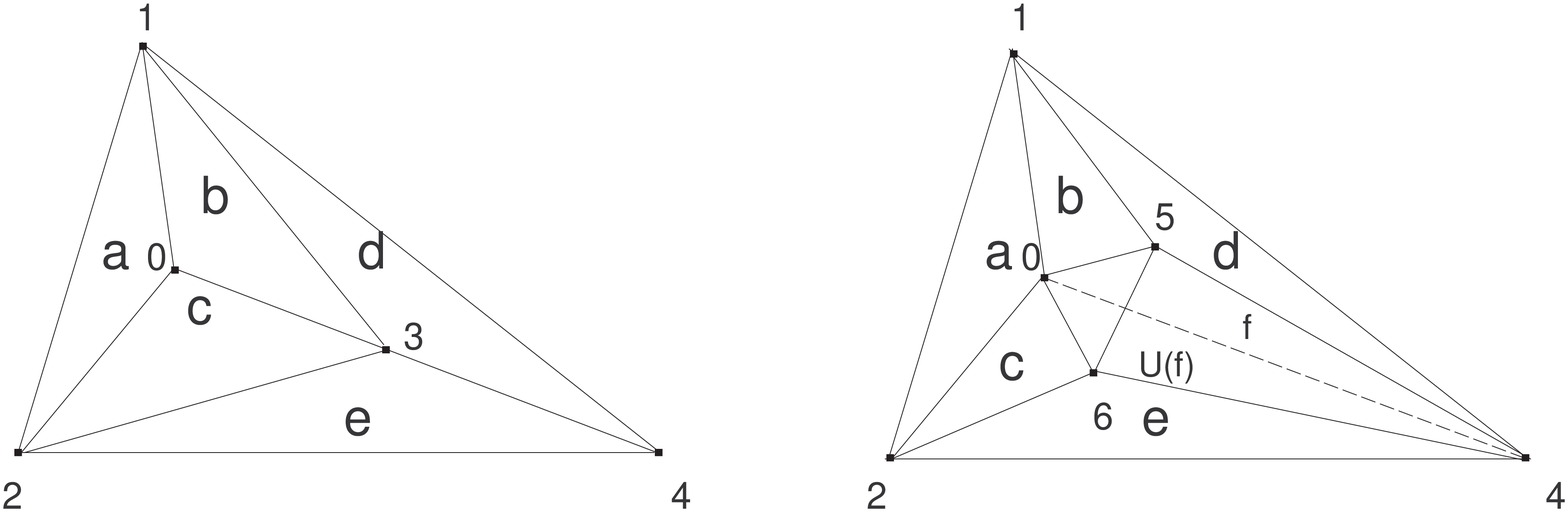,height=8cm,width=13cm}} 
\caption{$m(m(a,b,c),d,e)=m(a,m(b,d,f),m(c,U(f),e)$}
\label{fig6}
\end{figure}

In Figure \ref{fig6} we have two triangulation of the the 3-dimensional ball $B_3$. On the left hand side we have  two 3-simplexes
$[0,2,1,3]$ and $[3,2,1,4]$ which are  glued along the face $[2,1,3]$. On the right hand side we have four 3-simplexes $[0,5,6,4]$, $[5,0,1,4]$, $[6,2,0,4]$,
and $[0,2,1,4]$ which are glued along the faces $[0,5,4]$, $[6,0,4]$, $[0,1,4]$ and $[2,0,4]$ respectively. 

We consider $a\in \Gamma(\alpha,\beta^{-1})$, $b\in\Gamma(\beta,\gamma^{-1})$, $c\in \Gamma(\beta^{-1}\alpha,\gamma^{-1}\beta)$, $d\in \Gamma(\gamma, \delta^{-1})$, $e\in \Gamma(\gamma^{-1}\alpha,\delta^{-1}\gamma)$ and 
$f\in \Gamma(\gamma^{-1}\beta,\delta^{-1}\gamma)$

We can define two maps on $B_3$ such that  $[2,1,0]$ is mapped in $a$, $[0,1,3(5)]$ in $b$, $[2,0,3(6)]$ in $c$, $[3(5),1,4]$ in $d$, $[2,3(6),4]$ in $e$ and $[0,5,4]$ in $f$. Moreover we send $[5,6]$ in $m_0$ and for every point $x\in [5,6]$ we send $[x,0]$ in $\gamma^{-1}\beta$ and $[4,x]$ in $\delta^{-1}\gamma$. It follows that $[6,0,4]$ must go to $Q(f)$.

It is obvious  that up to homotopy the image of $[2,1,4]$ from the two maps is the same element in $\Gamma(\alpha,\delta^{-1})$. More exactly we have:
$$m(m(a,b,c),d,e)=m(a,m(b,d,f),m(c,Q(f),e)$$
It is easy to see that:
$$P(m(a,b,c))=m(P(b),P(c),P(a))$$
$$Q(m(b,a,f))=m(Q(b),Q(f),Q(a))$$
We want to prove an analog for formula \ref{s4}. The staring point is again Example \ref{dwe} where we take $\alpha$ to be the trivial 3-cocycle. If we  write the condition \ref{s4} for $a=(u,v,t)$ and $b=(x,y,z)\in k[G^{(3-1)}]$ we get that for every $k_1=(g_1,h_1,(h_1g_1)^{-1})$ there is an element  $k_2=(g_2,h_2,(h_2g_2)^{-1})$  such that the following equalities hold:
\begin{eqnarray*}
&&m((y,z,x),(u,v,t),(g_1,h_1,(h_1g_1)^{-1}))=(g_2,h_2,(h_2g_2)^{-1})\\
&&(g_1 ^{-1},h_1g_1,h_1 ^{-1})=m((x,y,z),(g_2^{-1},h_2g_2,h_2^{-1}),(u,v,t))
\end{eqnarray*}
This can be written as $m(P(b),a,k_1)=k_2$ and $Q(k_1)=m(b,Q(k_2),a)$. 
We combine the two equalities to get:
\begin{eqnarray*}
Q(k_1)=m(b,Q(m(P(b),a,k_1)),a)
\end{eqnarray*}
or after some changes of variables:
\begin{eqnarray}
f=m(m(f,a,b),P^2Q(a),PQ(b))
\label{s4gr}
\end{eqnarray}

One can see that equation \ref{s4gr} is  also true for every $f\in \Gamma(\alpha, \beta^{-1})$, $a\in \Gamma(\beta,\gamma^{-1})$ and $b\in 
\Gamma(\beta^{-1}\alpha,\gamma^{-1}\beta)$. 

We ready to give the formal definition of a $\Delta$-group. 

\begin{definition}
Let $G$ be a group. A $\Delta$-group  based at $G$ is a collection of sets 
$T= \{T(g,h)\}_{g,h\in G}$ together with operations:
\begin{eqnarray*}
&m:T(g,h^{-1})\times T(h,k^{-1})\times T(h^{-1}g,k^{-1}h)\to T(g,k^{-1})\\
&P:T(g,h)\to T(h,g^{-1}h^{-1})\\
&Q:T(g,h)\to T(g^{-1},hg)
\end{eqnarray*}
such that the following compatibilities hold:
\begin{eqnarray}
P^3=id, \; Q^2=id,\; P^2Q=QP \label{cm}
\end{eqnarray}
\begin{eqnarray}
P(m(a,b,c))=m(P(b),P(c),P(a))\label{pm}
\end{eqnarray}
\begin{eqnarray}
Q(m(a,b,c))=m(Q(a),Q(c),Q(b))\label{qm}
\end{eqnarray}
\begin{eqnarray}
m(m(a,b,c),d,e)=m(a,m(b,d,f),m(c,Q(f),e))\label{mm}
\end{eqnarray}
\begin{eqnarray}
m(m(f,a,b),P^2Q(a),PQ(b))=f\label{im}
\end{eqnarray}
\end{definition}
\begin{example} If $(A,+)$ is a commutative group, we define 
\begin{eqnarray*}
m(a,b,c)=a+b+c,\; P(a)=a,\; Q(a)=-a
\end{eqnarray*}\label{triv0}
then $T(\{1\},A)=\{T(1,1)=A\}$ becomes a $\Delta$-group based at the trivial group $\{1\}$.   
\end{example}
\begin{example}
If $G$ is a group set 
$T(g,h)=\{(g,h,(hg)^{-1})\}$ for all $g,h\in G$ and 
$T(G,0)=\{T(g,h)\}_{g,h\in G}$.  Then $T(G,0)$ is a $\Delta$-group based at $G$ with:
\begin{eqnarray*}
&m((g,h^{-1},g^{-1}h),(h,k^{-1},h^{-1}k),(h^{-1}g,k^{-1}h,g^{-1}k))=
(g,k^{-1},g^{-1}k)\\
&P((g,h,g^{-1}h^{-1}))= (h,g^{-1}h^{-1},g)\\
&Q((g,h,g^{-1}h^{-1}))= (g^{-1},hg, h^{-1})
\end{eqnarray*} \label{trivg}
We call $T(G,0)$ the trivial $\Delta$-group based at $G$.
\end{example}

\begin{example}
Let $G$ be a group and $(A,+)$ a $G$-module. Set    
$T(g,h)=\{(a,(g,h))\vert a\in A \}$ for all $g,h\in G$ and $T(G,A)=\{T(g,h)\}_{g,h\in G}$. Then $T(G, A)$ is a $\Delta$-group based at $G$ with:
\begin{eqnarray*}
&m(((a,(g,h^{-1})),(b,(h,k^{-1})),(c, (h^{-1}g,k^{-1}h)))=
(a+(g^{-1}h)b+c,(g,k^{-1}))\\
&P((a,(g,h^{-1})))= (ga,(h^{-1},g^{-1}h))\\
&Q((a,(g,h^{-1})))= (-(ga),(g^{-1},h^{-1}g))
\end{eqnarray*} \label{tga}
\end{example}

In the next section we shall give a more general example $T(G,A,\alpha)$ corresponding to special 3-cocycles $\alpha\in C^3(G,A)$, so we postpone the verification until then. Also, at that point, it will become clear how these examples where conceived.
 
We can now formulate the main result of this section:
\begin{theorem}
Let $M$ be a manifold with the property that $\pi_1(M)$ has no element of order two. Then $\Gamma(M)$ is a $\Delta$-group based at $\pi_1(M)$. Moreover if $f:M\to N$ is a map between two such manifolds then $f^*:\Gamma(M) \to\Gamma(N)$ is a morphism of $\Delta$-groups. 
\end{theorem}
\begin{proof}
It follows from the above construction.
\end{proof}
\begin{remark}
It is easy to see that $\Gamma(M)$ does not depend on  the set of paths ${\cal B}(M)$. If  ${\cal B}'(M)$ is another set of based curves, we can take a set of homotopies between pair of elements from ${\cal B}(M)$ and ${\cal B}'(M)$ and then using these homotopies we can construct an isomorphism between  $\Gamma'(M)$ and $\Gamma(M)$.
\end{remark}

\begin{remark} $\Gamma(M)$ can be defined for every path connected topological space.
\end{remark}
\begin{remark} Let $M$  be a manifold. If $\pi_1(M)=1$ then $\Gamma(M)=T(1,\pi_2(M))$ as in Example \ref{triv0}. If $\pi_2(M)=0$ then 
$\Gamma(M)=T(\pi_1(M),0)$ as in Example \ref{trivg}.
\label{pitriv}
\end{remark}

We end this section with an example of a $3$-algebra associated to a finite $\Delta$-group. This is similar with the construction of the group algebra $kG$  associated to a group $G$. 
\begin{example} Suppose that $T$ is finite a $\Delta$-group based at $G$. Define 
$A=\coprod_{g,h\in G}(\coprod_{x\in T(g,h)}kx)$. We extend $m$ and $P$ linearly  to the whole vector space $A$. Define $u:k\to A\otimes A$, 
$$u(1)=\sum_{g,h\in G}\frac{1}{\#T(g,h)}\sum_{x\in T(g,h)}x\otimes Q(x)$$
Straightforward computations show that $(A,m,u,P)$ is a strong 3-algebra.
\end{example}

\section{$\Gamma(M)\simeq T(\pi_1(M),\pi_2(M),\alpha)$}

In this section we study the structure of $\Gamma(M)$. We will prove that it is determined by the action of $\pi_1(M)$ on $\pi_2(M)$ and a certain 3-cocycle. 

In what follows $\pi_1(M)$ has a multiplicative operation, $\pi_2(M)$ has an additive operation, $g,h,...$  are elements from $\pi_1$ and  $a,b,...,f$ are elements of $\pi_2$. These conventions allows us to use $ga$ for the action of $\pi_1$ on $\pi_2$ without confusion. 
  
\begin{figure}[htbp] 
\vspace*{-8pt}
\centerline{\psfig{file=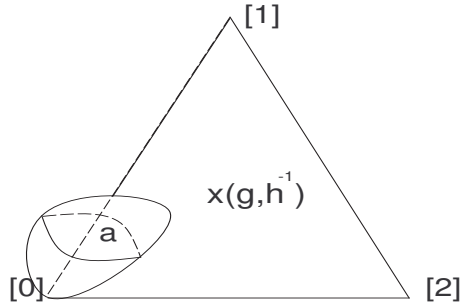,height=4cm,width=6cm}} 
\caption{$(a,(g,h^{-1}))$ a generic element from $\Gamma(g, h^{-1})$ }
\label{fig8}
\end{figure}

Let $M$ be a manifold such that there is no element of order 2 or 3 in $\pi_1(M)$. Fix an element $x(g,h^{-1})$ in each $\Gamma(g, h^{-1})$.  We may assume without lose of the generality that $P(x(g,h^{-1}))=x(h^{-1},g^{-1}h)$ and 
$Q(x(g,h^{-1}))=x(g^{-1},h^{-1}g)$. Notice that any other element $y\in \Gamma(g, h^{-1})$ differs from $x(g,h^{-1})$ by an element of $\pi_2(M)$; the only problem is where do we glue this bubble. By convention we assume that the element of $\pi_2$ is always at the $[0]$ corner of our two simplex (see Figure 
\ref{fig8}). This means that we can identify the set $\Gamma(M)(g,h^{-1})$ with $\pi_2(M)$. For convenience we denote such an element by $(a,(g,h^{-1}))$. 

If one wants to move the bubble from  the $[0]$ corner to the $[1]$ corner then one has to take in consideration the action of $\pi_1(M)$ on $\pi_2(M)$. 
Having this in mind, it is easy to see that:
$$P((a,(g,h^{-1})))=(ga,(h^{-1}, g^{-1}h))$$
$$Q((a,(g,h^{-1})))=(-ga,(g^{-1}, h^{-1}g))$$
\begin{figure}[htbp] 
\hspace*{-60pt}
\centerline{\psfig{file=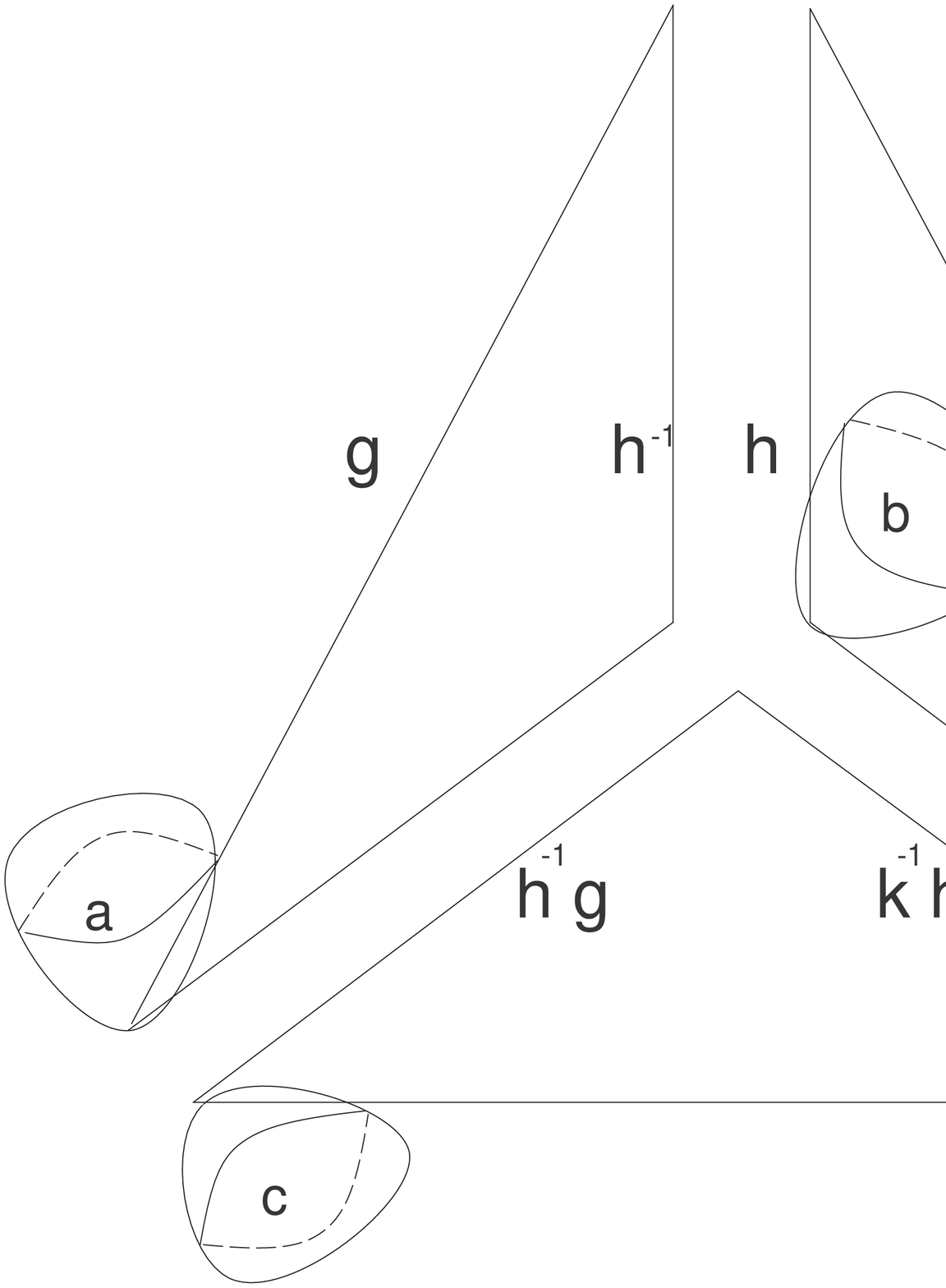,height=6cm,width=4.5cm}} 
\caption{$m((a,(g,h^{-1})),(b,(h,k^{-1})),(c, (h^{-1}g,k^{-1}h)))$}
\label{fig9}
\end{figure}

To find the multiplication, we look to Figure \ref{fig9}.  
There are two distinct problems. The first one is how to multiply the $x(g,h^{-1})$'s  among them and the second how to add the elements of $\pi_2(M)$. Fortunately the two problems are independent.  

First take  $x(g,h^{-1}),x(h,k^{-1}),x(h^{-1}g,k^{-1}h)\in \Gamma(M)$. The product of these three elements belong to $\Gamma(M)(g,k^{-1})$ so it is just $x(g,k^{-1})$ plus an element $y(g,h,k)\in \pi_2(M)$ which is glued in the $[0]$ corner. We define $\alpha:G\times G\times G\to A$ by $\alpha(g^{-1}h,h^{-1}k,k^{-1})=y(g,h,k)$. To be more precise we have:
\begin{eqnarray}
\alpha(g,h,k)=y((ghk)^{-1},(hk)^{-1},k^{-1}) \label{3cm}
\end{eqnarray}
We will prove later that this is a 3-cocycle.

Now let  $(a,(g,h^{-1})),(b,(h,k^{-1})),(c, (h^{-1}g,k^{-1}h))\in \Gamma(M)$. The $\pi_2(M)$ component of the first and third element are already in the $[0]$ corner. However for the second element we have to move $b$ along $g^{-1}h$. This means that the $\pi_2(M)$ contribution to the product is: $a+(g^{-1}h)b+c$. To conclude we have that in $\Gamma(M)$ the multiplication is defined by:
$$m((a,(g,h^{-1})),(b,(h,k^{-1})),(c, (h^{-1}g,k^{-1}h)))=
(a+(g^{-1}h)b+c+\alpha(g^{-1}h,h^{-1}k,k^{-1}),(g,k^{-1}))$$
It is clear now how we constructed Example \ref{tga} and why Remark \ref{pitriv} is true. 

More general, suppose that $G$ is a group, $A$ a $G$-module and $\alpha:G\times G\times G\to A$. We want to  see under what conditions the above maps define a $\Delta$-group $T(G,A,\alpha)$.
 
First we check that $P^3=id$.   
\begin{eqnarray*}P^3((a,(g,h^{-1})))&=&P^2((ga,(h^{-1},g^{-1}h)))\\
&=&P((h^{-1}ga,(g^{-1}h,g)))\\
&=&(g^{-1}hh^{-1}ga,(g,h^{-1}))\\
&=&(a,(g,h^{-1}))
\end{eqnarray*}  
Similarly one has that $Q^2=id$ and $Q^2P=PQ$.

Let's check (\ref{pm}).
\begin{eqnarray*}
&&P(m((a,(g,h^{-1})),(b,(h,k^{-1})),(c, (h^{-1}g,k^{-1}h))))=\\
&=&P((a+(g^{-1}h)b+c+\alpha(g^{-1}h,h^{-1}k,k^{-1},(g,k^{-1})))\\
&=&(g(a+(g^{-1}h)b+c+\alpha(g^{-1}h,h^{-1}k,k^{-1})),(k^{-1},g^{-1}k))\\
&=&(ga+hb+gc+g\alpha(g^{-1}h,h^{-1}k,k^{-1}),(k^{-1},g^{-1}k))
\end{eqnarray*}
\begin{eqnarray*}
&&m(P(b,(h,k^{-1})),P(c, (h^{-1}g,k^{-1}h)),P(a,(g,h^{-1})))=\\
&=&m((hb,(k^{-1}, h^{-1}k)),(h^{-1}gc, (k^{-1}h,g^{-1}k)),(ga,(h^{-1},g^{-1}h)))\\
&=&(hb+k(h^{-1}k)^{-1}h^{-1}gc+ga+\alpha(h,h^{-1}g,g^{-1}k),(k^{-1},g^{-1}k))\\
&=&(hb+gc+ga+\alpha(h,h^{-1}g,g^{-1}k),(k^{-1},g^{-1}k))
\end{eqnarray*}

And so, in order to have (\ref{pm}), $\alpha$ must satisfy the following identity:
\begin{eqnarray}
g\alpha(g^{-1}h,h^{-1}k,k^{-1})=\alpha(h,h^{-1}g,g^{-1}k) \label{c1}
\end{eqnarray} 
A similar computation shows that (\ref{qm}) is equivalent with:
\begin{eqnarray}
g\alpha(g^{-1}h,h^{-1}k,k^{-1})=-\alpha(h,h^{-1}k,k^{-1}g) \label{c2}
\end{eqnarray} 
Let's look to (\ref{mm}).
\begin{eqnarray*}
&&m(m((a,(g,h^{-1})),(b,(h,k^{-1})),
(c, (h^{-1}g,k^{-1}h))),(d,(k,l^{-1})),(e,(k^{-1}g,l^{-1}k)))=\\
&=&m((a+g^{-1}hb+c+\alpha(g^{-1}h,h^{-1}k,k^{-1}),(g,k^{-1})),(d,(k,l^{-1})), (e,(k^{-1}g,l^{-1}k)))\\
&=&(a+g^{-1}hb+c+\alpha(g^{-1}h,h^{-1}k,k^{-1})+g^{-1}kd+e+\alpha(g^{-1}k,
k^{-1}l,l^{-1}),(g,l^{-1}))
\end{eqnarray*}
\begin{eqnarray*}
&&m((a,(g,h^{-1})),m((b,(h,k^{-1})),(d,(k,l^{-1})),(f,(k^{-1}h,l^{-1}k))),
m((c, (h^{-1}g,k^{-1}h)),\\
&&Q((f,(k^{-1}h,l^{-1}k))),(e,(k^{-1}g,l^{-1}k))))=\\
&=&m((a,(g,h^{-1})),(b+h^{-1}kd+f+\alpha(h^{-1}k,k^{-1}l,l^{-1}),(h,l^{-1})),
(c-(g^{-1}h)f+e+\\
&&\alpha(g^{-1}k,k^{-1}l,l^{-1}h), (h^{-1}g,l^{-1}h)))\\
&=&(a+(g^{-1}h)(b+h^{-1}kd+f+\alpha(h^{-1}k,k^{-1}l,l^{-1}))+c-(g^{-1}h)f+e+\\
&&\alpha(g^{-1}k,k^{-1}l,l^{-1}h)
+\alpha(g^{-1}h,h^{-1}l,l^{-1}),(g,l^{-1}))
\end{eqnarray*}
And so (\ref{mm}) is equivalent with:
\vspace{10pt}

\lefteqn{\alpha(g^{-1}h,h^{-1}k,k^{-1})+\alpha(g^{-1}k,
k^{-1}l,l^{-1})=}
\vspace*{-20pt}
\begin{eqnarray}
=g^{-1}h\alpha(h^{-1}k,k^{-1}l,l^{-1}))+\alpha(g^{-1}k,k^{-1}l,l^{-1}h)+\alpha(g^{-1}h,h^{-1}l,l^{-1}) \label{c3}
\end{eqnarray}
Finally one can prove that (\ref{im}) is equivalent with:
\begin{eqnarray}
\alpha(g^{-1}h,h^{-1}k,k^{-1})=-\alpha(g^{-1}k,k^{-1}h,h^{-1})
\label{c4}
\end{eqnarray}
Making an appropriate change of variables (\ref{c1}), (\ref{c2}) and (\ref{c4}) can be written as:
\begin{eqnarray}
\alpha(x,y,z)=xy\alpha(y^{-1},yz,(xyz)^{-1})=-\alpha(x,yz,z^{-1})=-\alpha(xy,y^{-1},yz) \label{d1}
\end{eqnarray}     
while (\ref{c3}) becomes:
\begin{eqnarray}
\alpha(x,y,zt)+\alpha(xy,z,t)=x\alpha(y,z,t)+\alpha(xy,z,(yz)^{-1})+\alpha(x,yz,t)\label{d2}
\end{eqnarray}
Using (\ref{d1}) twice we get:
\begin{eqnarray*}
\alpha(x,y,z)=\alpha(xy,z,(yz)^{-1})
\end{eqnarray*}    
and so $\alpha$ is a 3-cocycle. To conclude we have proved:

\begin{proposition} $T(G,A,\alpha)$ is a $\Delta$-group if and only if $\alpha\in Z^3(G,A)$ and (\ref{d1}) holds.
\end{proposition}  
\hfill $ \square $ 
\begin{remark} For $G$, $A$ and $\alpha$ as above we have an exact sequence of $\Delta$-groups:
\begin{eqnarray*}
0\to T(1,A)\to T(G,A,\alpha)\to T(G,0)\to 1
\end{eqnarray*} 
\end{remark}
     
Let's take $f:T(G,A,\alpha)\to  T(G,A,\beta)$ which is compatible with the  short exact sequence. This means that:  
$$f((a,(g,h^{-1})))=(a+\sigma(g^{-1}h, h^{-1}),(g,h^{-1}))$$
where $\sigma:G\times G\to A$. Because $f$ is a morphism of $\Delta$-groups, $f$ must be compatible with $P$, $Q$ and $m$. This yields the following conditions:
\begin{eqnarray}
\sigma(g,h)=g\sigma(h,(gh)^{-1})=-(gh)\sigma(h^{-1},g^{-1})
\label{e1}
\end{eqnarray}
\begin{eqnarray}
\sigma(gh,k)+\alpha(g,h,k)=\sigma(h,hk)+g\sigma(h,k)+\sigma(gh,h^{-1})+\beta(g,h,k)
\label{e2}
\end{eqnarray}

From (\ref{e1}) it follows that 
\begin{eqnarray*}
\sigma(g,h)=-\sigma(gh,h^{-1})
\end{eqnarray*} 
and so (\ref{e2}) says that $\alpha$ and $\beta$ are equivalent 3-cocycles. To be a little more precise we have to construct the first few terms of the  ``symmetric cohomology".

Let $G$ be a group, $A$ a $G$-module and $C^n(G,A)=\{\alpha:G^n\to A\}$. For $n=1,2$ or $3$ we have an action of the symmetric group 
$\Sigma_{n+1}$ on $C^n(G, A)$:

If $\phi\in C^1(G,A)$
\begin{eqnarray*}
((1,2)\phi)(g)=-g\phi(g^{-1})
\end{eqnarray*}
if $\sigma\in C^2(G,A)$
\begin{eqnarray*}
((1,2)\sigma)(x,y)=-x\sigma(x^{-1},xy)\\
((2,3)\sigma)(x,y)=-\sigma(xy,y^{-1})
\end{eqnarray*}
if $\alpha\in C^3(G, A)$
\begin{eqnarray*}
((1,2)\alpha)(x,y,z)=-x\alpha(x^{-1},xy,z)\\
((2,3)\alpha)(x,y,z)=-\alpha(xy,y^{-1},yz)\\
((3,4)\alpha)(x,y,z)=-\alpha(x,yz,z^{-1})
\end{eqnarray*}

In each case we define $CS^n(G, A)=(C^n(G,A))^{\Sigma_{n+1}}$. It is easy to see that $CS^n(G,A)$ is a subcomplex of the usual cohomology complex. We define: 
$$HS^n(G,A)=\frac{ZS^n(G, A)}{BS^n(G,A)}$$
With these notations we have the following result:
\begin{theorem} $f:T(G,A,\alpha)\to  T(G,A,\beta)$ is an isomorphism if and only if $[\alpha]=[\beta]$ in $HS^3(G,A)$.
\end{theorem}
\begin{proof} Straightforward.
\end{proof}
\begin{corollary}
The element $[\alpha]\in HS^3(\pi_1(M),\pi_2(M))$ defined by (\ref{3cm}) is an invariant of the space $M$. 
\end{corollary}

\section{Symmetric cohomology for groups}

Let $G$ be a group and $A$ a $G$ module. As usual we define $C^n(G,A)=\{\sigma : G^n\to A\}$, $\partial_n: C^n(G,A)\to C^{n+1}(G,A)$
\begin{eqnarray*}
\partial_n(\sigma)(g_1,...,g_{n+1})=g_1\sigma(g_2,...,g_{n+1})
-\sigma(g_1g_2,g_3,...,g_{n+1})+...+\\
+(-1)^n\sigma(g_1,...,g_ng_{n+1})+(-1)^{n+1}\sigma(g_1,...,g_n)
\end{eqnarray*} 
 Define $d^j:C^n(G,A)\to C^{n+1}(G,A)$ by 
\begin{eqnarray*}
&&d^0(\sigma)(g_1,...,g_{n+1})=g_1\sigma(g_2,...,g_{n+1})\\
&&...\\
&&d^j(\sigma)(g_1,...,g_{n+1})=\sigma(g_1,...g_jg_{j+1},...,g_{n+1})\\
&&...\\
&&d^{n+1}(\sigma)(g_1,...,g_{n+1})=\sigma(g_1,...,g_n)
\end{eqnarray*}
Let's notice that $\partial_n(\sigma)=\sum_0^{n+1}(-1)^jd^j$.
It is well known that in this way we obtain a complex and it's homology groups are denoted with $H^n(G,A)$.
We give here an action of $\Sigma_{n+1}$ on $C^n(G,A)$ (for every $n$) and  prove that it is compatible with the differential. This allows to develop the whole theory of the symmetric cohomology. 

It's enough to say what is the action of the transposition $(i,i+1)$ for $1\leq i\leq n$. For $\sigma \in C^n(G,A)$ we define:
\begin{eqnarray*}
&&((1,2)\sigma)(g_1,g_2,g_3,...,g_n)=-g_1\sigma((g_1)^{-1},g_1g_2,g_3,...,g_n)\\
&&((2,3)\sigma)(g_1,g_2,g_3,...,g_n)=-\sigma(g_1g_2,(g_2)^{-1},g_2g_3,g_4,...,
g_n)\\
&&...\\
&&((n-1,n)\sigma)(g_1,g_2,g_3,...,g_n)=-\sigma(g_1,g_2,..., g_{n-2}g_{n-1},(g_{n-1})^{-1},g_{n-1}g_n)\\
&&((n,n+1)\sigma)(g_1,g_2,g_3,...,g_n)=-\sigma(g_1,g_2,g_3,...g_{n-1}g_n,
(g_n)^{-1})
\end{eqnarray*}

\begin{proposition} The above formulas define an action of $\Sigma_{n+1}$ 
on $C^n(G,A)$ which is compatible with the differential $\partial$.
\label{scg}
\end{proposition}
\begin{proof} Let's see that we have an action of $\Sigma_{n+1}$ on $C^n(G,A)$.  First we check that the square of the action of the transposition $(i,i+1)$ is the identity.
\begin{eqnarray*} 
&&((i,i+1)((i,i+1)\sigma)))(g_1,g_2,... ,g_n)\\
&=&-((i,i+1)\sigma)(g_1,...,g_{i-1}g_i,g_i^{-1},g_ig_{i+1},...,g_n)\\
&=&-(-(\sigma)(g_1,..., g_{i-1}g_ig_i^{-1},(g_i^{-1})^{-1},g_i^{-1}
g_ig_{i+1},...,g_n))\\
&=&\sigma(g_1,...,g_n)
\end{eqnarray*}
For the braid relation we have:
\begin{eqnarray*} 
&&((i,i+1)((i+1,i+2)((i,i+1)\sigma)))(g_1,g_2,... ,g_n)\\
&=&-((i+1,i+2)((i,i+1)\sigma))(g_1,...,g_{i-1}g_i,g_i^{-1},g_ig_{i+1},...,g_n)\\
&=&((i,i+1)\sigma)(g_1,...,g_{i-1}g_i,g_i^{-1}g_ig_{i+1},(g_ig_{i+1})^{-1},g_ig_{i+1}g_{i+2},...,g_n)\\
&=&-\sigma(g_1,...,g_{i-1}g_ig_{i+1},g_{i+1}^{-1},g_{i+1}(g_ig_{i+1})^{-1},g_ig_{i+1}g_{i+2},...,g_n)\\
&=&-\sigma(g_1,...,g_{i-1}g_ig_{i+1},g_{i+1}^{-1},g_i^{-1},g_ig_{i+1}g_{i+2},...,g_n)
\end{eqnarray*}
and similarly 
\begin{eqnarray*} 
&&((i+1,i+2)((i,i+1)((i+1,i+2)\sigma)))(g_1,g_2,... ,g_n)\\
&=&-\sigma(g_1,...,g_{i-1}g_ig_{i+1},g_{i+1}^{-1},g_i^{-1},g_ig_{i+1}g_{i+2},...,g_n)
\end{eqnarray*}
So $(i,i+1)((i+1,i+2)((i,i+1)\sigma))=(i+1,i+2)((i,i+1)((i+1,i+2)\sigma))$.
All the other relations are easy to check.

We also want to prove that this action is compatible with $\partial$. More exactly, if $\sigma\in C^n(G,A)$ is invariant under the action of $\Sigma_{n+1}$ then $\partial(\sigma)$ is invariant under the action of $\Sigma_{n+2}$. 
We can check that 
\begin{eqnarray*}
&&(i,i+1)(d^j(\sigma))=d^j((i,i+1)\sigma)\;  {\rm if} \; i\leq j\\ &&(i,i+1)(d^j(\sigma))=d^j((i-1,i)\sigma) \; {\rm if} \; j+2\leq i\\
&&(i,i+1)(d^{i-1}(\sigma))=-d^i(\sigma)\\
&&(i,i+1)(d^i(\sigma))=-d^{i-1}(\sigma)
\end{eqnarray*}
Now if we take $\sigma$ which is invariant by $\Sigma_{n+1}$ and use the fact that $\partial(\sigma)=\sum_0^{n+1}(-1)^jd^j$ we get  $(i,i+1)(\partial(\sigma))=\partial(\sigma)$, which finish the proof. 
\end{proof}
\begin{definition}
The subcomplex obtained in Proposition \label{sdc} is denoted by $CS^n(G,A)$.  Its homology, $HS^n(G,A)$, is called the symmetric cohomology of $G$ with coefficients in $A$.
\end{definition}
\begin{remark}
There is a natural map from $HS^n(G,A)$ to $H^n(G,A)$. When $n=1$ or $n=2$ it is easy to check that the map is injective. However for $n\geq3$ it is not clear if this fact is still true. 
\end{remark}

%

\end{document}